\documentclass[12pt, reqno]{amsart}
\usepackage{amsmath, amsthm, amscd, amsfonts, amssymb, graphicx, color}
\usepackage{setspace}
\usepackage{mathrsfs}
\usepackage{multicol}
\usepackage{tikz-cd}
\usepackage[bookmarksnumbered, colorlinks, plainpages]{hyperref}
\hypersetup{colorlinks=true,linkcolor=red, anchorcolor=green, citecolor=cyan, urlcolor=red, filecolor=magenta, pdftoolbar=true}

\newtheorem{theorem}{Theorem}[section]
\newtheorem{lemma}[theorem]{Lemma}

\newtheorem{corollary}[theorem]{Corollary}
\theoremstyle{definition}

\theoremstyle{remark}

\numberwithin{equation}{section}

\newcommand{\NN}{\mathbb{N}}

\newcommand{\CC}{\mathbb {C}}

\newcommand{\R}{\mathbb{R}}

\newcommand{\D}{\mathbb{D}}
\begin{document}
\setcounter{page}{1}
\title[Convex-cyclic weighted composition operators and their adjoints   ]{Convex-cyclic weighted composition operators and their adjoints}
\author[Tesfa  Mengestie]{Tesfa  Mengestie }
\address{ Mathematics Section \\
Western Norway University of Applied Sciences\\
Klingenbergvegen 8, N-5414 Stord, Norway}
\email{Tesfa.Mengestie@hvl.no}


 \begin{abstract}
 We characterize the convex-cyclic weighted composition operators $W_{(u,\psi)}$ and their adjoints  on the  Fock space  in terms of   the derivative powers  of $ \psi$ and the location of  the eigenvalues of the operators on the complex plane. Such a description is also equivalent to identifying  the   operators or their adjoints for which their  invariant closed convex sets are all invariant subspaces.    We further show that the space supports no supercyclic weighted composition operators with  respect to the pointwise convergence topology and hence with the weak and strong  topologies and answers  a question raised by   T. Carrol and C. Gilmore in \cite{CC}.\\


\end{abstract}
\maketitle
\section{Introduction}
The study of the weighted composition operator $ W_{(u,\psi)}:  f\mapsto  u\cdot  f( \psi) $  with symbol $\psi$ and multiplier $u$ acting on various spaces of holomorphic functions traces back to  works  related to isometries on the Hardy spaces  \cite{For,Hof}  and commutants of Toeplitz operators \cite{CC1,CC2}. Since then the  operator has attracted much research interest and there exists now   rich body of literatures  dealing with many  of its properties in various settings; see for example \cite{BM2, CO,Tle,TM4} and the references therein.

In this note,  we  are interested in    linear dynamical properties of the operators and their adjoints   on the  Fock space $\mathcal{F}_2$ which   consists of  square integrable analytic functions in $\CC$ with respect to the Gaussian measure $d\mu(z)= \frac{1}{\pi}e^{-|z|^2}  dA(z)$ where $dA$ is the Lebesgue measure in  $\CC$.  The space  is a reproducing kernel Hilbert space  endowed with the inner product
\begin{align*}
\langle f, g\rangle=   \int_{\CC} f(z)\overline{g(z)}d\mu(z),
\end{align*}norm  $\|f\|_2:=\sqrt{ \langle f, f\rangle} $ and  kernel function $K_w(z)= e^{\langle z, w\rangle}$.

A bounded linear operator $T$ on a separable  Banach space $\mathcal{H}$ is said to be cyclic  if there exists a vector $f$ in $\mathcal{H}$ for which the   span of the orbit  \begin{align*}\text{Orb}(T,f)=\big\{ f,\  Tf,\  T^2f,\  T^3 f, \  ...\big\}\end{align*}  is dense in $\mathcal{H}$. Such an $f$ is called a cyclic vector  for $T$.  The operator is hypercyclic if  the   orbit it self  is dense, and  supercyclic with vector $f$ if  the projective orbit,
 \begin{align*}\text{ Projorb}(T,f)= \big\{ \lambda T^nf, \ \ \lambda \in \CC, \ n= 0, 1, 2, ...\big\},\end{align*}
 is dense.
 These dynamical properties of  $T$ depend on the behaviour of its iterates
 \begin{align*}
 T^n= \underbrace{T \circ T\circ T \circ...\circ T}_{n - \text{times}}\end{align*}
 We may note that identifying cyclic and  hypercyclic operators have been  a subject of high interest partly because they play central rolls in the study of other operators. More specifically, it is known
that  every bounded linear operator
on an infinite dimensional  complex separable  Hilbert space is the sum of two hypercyclic
operators \cite[p. 50]{BM}. This result holds true with the summands being cyclic operators  as well \cite{WU}.
\subsection{Cyclic weighted composition operators and their adjoints}
The bounded and compact properties of   weighted composition operators on $\mathcal{F}_2$ are identified in \cite{Tle,TM4},  and    $W_{(u,\psi)}$ is bounded if and only if $u \in \mathcal{F}_2$  and
 \begin{align}
 \label{bounded}
 \sup_{z\in \CC} |u(z)|e^{\frac{1}{2}(|\psi(z)|^2-|z|^2)} <\infty.
 \end{align}
  Furthermore, it  was   proved that condition \eqref{bounded} implies $
 \psi(z)= az+b, |a|\leq 1$ and when $a\neq 0$ the operator norm is estimated by
  \begin{align}
 \label{normest}
 \sup_{z\in \CC} |u(z)|e^{\frac{1}{2}(|az+b|^2-|z|^2)} \leq \|W_{(u,\psi)}\|\leq |a|^{-1} \sup_{z\in \CC} |u(z)|e^{\frac{1}{2}(|az+b|^2-|z|^2)}.
 \end{align}
 When $|a|=1$, the multiplier function  has the special  form
 $
 u=u(0)K_{-\overline{a}b}
$ and the relation in \eqref{normest} simplifies to
  \begin{align}\label{simplenorm}
 \|W_{(u,\psi)}\|= |u(0)|e^{\frac{|b|^2}{2}}.
 \end{align}
  In \cite{TMW},  we reported that there exists no supercyclic  composition operator on Fock spaces. On the other hand,  the orbit of any vector $f$ under  $W_{(u,\psi)}$  has elements of the form
\begin{align}
\label{interplay}
W_{(u,\psi)}^n f= f(\psi^n) \prod_{j=0}^{n-1} u(\psi^j)
\end{align} for all nonnegative integers $n$ and $\psi^0$ is the identity map. The relation in \eqref{interplay} shows that the power of the  weighted composition operators is another weighted composition operator with symbol  $(u_n, \psi^n)$ where
\begin{align}
\label{prod}
 u_n= \prod_{j=0}^{n-1} u(\psi^j).
\end{align}
In \cite{TMW5},  we continued the work in \cite{TMW} and investigated whether the relation in \eqref{interplay} enduces  an interplay between the symbol $\psi $ and the multiplier  $u$ and result in  supercyclic  weighted composition operators.  It turns out that  such an interplay fails to make any projective orbit  dense  enough in  $\mathcal{F}_2$. See also \cite{CC} for a different approach.  Recently,  the author  \cite{TM6}   proved that the Fock space support no supercyclic  adjoint weighted  composition operators either. Recall that, the adjoint $W_{(u,\psi)}^*$ of a bounded  weighted composition operator $W_{(u,\psi)}$ on $\mathcal{F}_2$ is the operator which satisfies the relation
 \begin{align*}
 \langle W_{(u,\psi)} f, g\rangle= \langle f, W_{(u,\psi)}^*g\rangle
 \end{align*} for all $f, g \in \mathcal{F}_2$.  We note that the  adjoint of a weighted composition operator on $ \mathcal{F}_2$ is not necessarily a weighted composition operator.  L. Zhao and C. Pang \cite{LHP} proved that for  pairs of entire functions $u_1, \psi_1$ and $u_2, \psi_2$ from $\mathcal{F}_2$,  the relation $W_{(u_1, \psi_1)}^*= W_{(u_2, \psi_2)}$  holds if and only if
 \begin{align*}
 \psi_1(z)= az+b, u_1(z)= dK_c(z), \psi_2(z)= \overline{a}z+c, \text{and} \ \ u_2(z)= \overline{d}K_c(z)
 \end{align*} where $a, b, c$ and $d$ are constants such that  $d\neq 0$ and either $|a|<1$  or $|a|=1$ and $c+\overline{a}b=0$. Thus, an operator and its adjoint can have quite different dynamical structures; see for example \cite[p.26]{BM}  about the hypercyclic structure of  the multiplication operator on Hardy spaces.

 Having observed the absence of   supercyclic weighted
 composition operators and their adjoints on the Fock space, we considered the cyclicity problem in  \cite{TM6} and proved   the following interesting result.
 \begin{theorem} \label{thm1}
   Let $u$ and $\psi$ be entire functions on $\CC$ such that $W_{(u,\psi)}$  is bounded on  $\mathcal{F}_2$ and hence  $\psi(z)=az+b, \ |a|\leq 1$.     Then the following statements are equivalent.
   \begin{enumerate}
\item  $W_{(u,\psi)}$
  is cyclic on $\mathcal{F}_2$;
      \item  $u$ is non-vanishing  and   $a^k\neq a$ for all positive integer $ k\geq2$;
  \item  $W_{(u,\psi)}^*$
  is cyclic on $\mathcal{F}_2$.
        \end{enumerate}
   \end{theorem}
In this note we plan  to  answer two other  basic  questions related to the dynamics of the operators. The first  deals with identifying  weighted composition operators which admit convex-cyclicity property in the space. As will be seen latter, convex-cyclicity is stronger than the cyclicity property and may require conditions stronger than  part (ii) Theorem~\ref{thm1}. The  second takes up the question whether  weakening the topology of the space results in  weakly supercyclic weighted composition operators.
  \subsection{Convex-cyclic  weighted composition operators and their adjoints } \ \ \\
    Another dynamical  concept  related to the iterates of an operator is convex-cyclicity.    A bounded operator $T$ on a Banach space $\mathcal{H}$ is said to be convex-cyclic if there is a vector $f$  in $  \mathcal{H}$  such that the convex hull of  $ \text{Orb}(T,f)$  is dense in
$\mathcal{H}$. Recall that  the convex hull of a set is the set of all convex combinations of its elements, that is, all finite
linear combinations of its elements  where the coefficients are
non-negative and their sum is one.  Thus,  studying the convex-cyclicity of an operator   requires a good understanding of  the convex combinations
of  its  powers. The notion of convex-cyclicity  is a relatively young subject of study which was introduced by Rezaei \cite{HR} in 2013. A few more studies  were made recently  in \cite{TAN,NFP,LMP}.

  Since the convex hull of a set lies between the set and its linear span,  every hypercyclic operator
is convex-cyclic  while  every convex-cyclic operator is cyclic, and obviously, every convex-cyclic vector for  an operator  is a cyclic vector.  On a finite dimensional Banach space, the notions of cyclicity and convex-cyclicity are equivalent \cite[Theorem 1.1]{NFP}. Suppercyclicity is another structure that stands between hypercyclic and cyclic operators and one may  wonder its location in reference to convex-cyclicity.  By \cite[Proposition 3.2]{HR}, the norm of every convex-cyclic operator is bigger than one.  Hence,  if a bounded operator $T$ is supercyclic, then for positive parameters $\alpha$, all the operators  $1/(\alpha+\|T\|) T$ are supercyclic but not convex-cyclic.  On the other hand, as shown in Theorem~\ref{thm22}, there exists convex-cyclic weighted composition operators on  $\mathcal{F}_2$ which are not supercyclic; see the diagram in the last section for a good illustration of the relations among the various forms of cyclicities.

Then, it naturally follows to ask  which cyclic weighted composition operators are convex-cyclic on $\mathcal{F}_2$. This note aims to answer  this question. Our first result, Theorem~\ref{thm22}, completely  describes  the convex-cyclic weighted composition operators and their adjoints   in  terms of  cyclicity and  the   location of the  eigenvalues of the operators  on the complex plane. The second main result, Theorem~\ref{supercyclic}, shows that  the space fails to support supercyclic weighted composition operators  with respect to  the pointwise topology and hence with the weak topology.
 \begin{theorem} \label{thm22}
   Let $u$ and $\psi$  be analytic maps on $\CC$  such that $W_{(u,\psi)}$  is bounded on  $\mathcal{F}_2$, and hence $\psi(z)=az+b, \ |a|\leq 1$.   Then the following statements are equivalent.
   \begin{enumerate}
     \item  $W_{(u,\psi)}$  is convex-cyclic on $\mathcal{F}_2$;
     \item  $W_{(u,\psi)}^*$  is convex-cyclic on $\mathcal{F}_2$;
     \item $W_{(u,\psi)}$
 has the property that all of its invariant closed convex-sets are invariant subspaces;
 \item $W_{(u,\psi)}^*$
 has the property that all of its invariant closed convex-sets are invariant subspaces;
\item  $W_{(u,\psi)}$
  is cyclic on $\mathcal{F}_2$,   $|a|=1$, $|u(z_0)|>1$, and $\Im\big(u(z_0) a^m\big)\neq 0$ for all $m\in\NN_0$ where $\Im(z)$ refers to the imaginary part of a complex number $z$ and $z_0: = b/(1-a)$.
        \end{enumerate}
   \end{theorem}
   The theorem describes the convex-cyclic weighted composition operators and their adjoints by simple to  check   requirements. The condition $ |a|=1$ restricts that a non-normal weighted composition operator can not be convex-cyclic on $\mathcal{F}_2$, while the condition $\Im\big(u(z_0) a^m\big)\neq 0$ requires  all the eigenvalues of the operator and its adjoint not to be located on the real line.

   The proof of the theorem follows from Lemma~\ref{lem22}, Lemma~\ref{lem33} and Lemma~\ref{lem44} where their proofs are mainly  based on the relation between the dynamical and
    spectral  properties of the
   operators. The relationship between the density of orbits and the spectral properties
of an operator plays a vital roll in the study of dynamical structures of operators. Thus, identifying the location of the  eigenvalues of the operator has been  a fundamental tool for studying convex-cyclic operators.   We plan to use this tool to prove our main results.
  \begin{lemma} \label{lem22}
   Let $u$ and $\psi$  be analytic maps on $\CC$  such that $W_{(u,\psi)}$  is bounded on  $\mathcal{F}_2$, and hence $\psi(z)=az+b, \ |a|\leq 1$.   Then $W_{(u,\psi)}$  is convex-cyclic if and only if  the following holds.
     \begin{enumerate}
\item  $W_{(u,\psi)}$
  is cyclic on $\mathcal{F}_2$;
      \item  $|a|=1$, $|u(z_0)|>1$ and $\Im\big(u(z_0) a^m\big)\neq 0$ for all $m\in \NN_0$.
        \end{enumerate}
   \end{lemma}
 \emph{ Proof}. Let us first proof the sufficiency.  By Theorem~\ref{thm1}, the cyclicity condition implies that $a^m \neq a$ for all $m\geq 2$.  By \cite[Lemma~2]{TM6},  the numbers  $a^m u(z_0)$ constitutes a sequence of distinct eigenvalues for the operator  with corresponding eigenvectors
 \begin{align}
 \label{egv}
 f_m(z)= (z-z_0)^m e^{\frac{a\overline{b}}{a-1}z}
 \end{align} for all $m\in \NN_0$. Furthermore, since the sequence of the polynomials  is dense in $\mathcal{F}_2$  and $e^{\frac{a\overline{b}}{a-1}z}$ is a non-vanishing function in  $\mathcal{F}_2$, by a result of K. H. Izuchi  \cite{IHH}, the sequence of the eigenvectors $(f_m), m\in \NN_0$ is also dense in  $\mathcal{F}_2$. This along with condition (ii) of the theorem and Theorem~6.2 in \cite{TAN} ensure  that the operator is convex-cyclic, and  indeed has a dense set  of convex-cyclic vectors.

Conversely, suppose now that the operator is convex-cyclic. Since every convex-cyclic operator is cyclic, condition (i) follows and hence  $a\neq0$ and $a\neq1$.  Then, a  simple modification, like  changing the corresponding kernel function, of the arguments in the proof of Lemma~3 of \cite{GG}  gives that the adjoint operator $W_{(u,\psi)}^*$  has  the following  set of eigenvalues:
\begin{align*}
\big\{\overline{u(z_0)a^m}: m\in \NN_0\big\}.
\end{align*}
By \cite[Proposition~3.3]{HR}, it follows that $\overline{u({z_0})a^m}\in \CC\backslash \big(\overline{\D}\cup \R\big)$. Therefore,
$\Im\big(u(z_0) a^m\big)\neq 0$ and $|a^m u(z_0)|= |a|^m |u(z_0)|>1$ for all $m\in \NN_0$ where the later holds only when $|a|=1$ and $|u(z_0)|>1$, and completes the proof.
\begin{lemma} \label{lem33}
   Let $u$ and $\psi$  be analytic maps on $\CC$  such that $W_{(u,\psi)}$  is bounded on  $\mathcal{F}_2$, and hence $\psi(z)=az+b, \ |a|\leq 1$.   Then $W_{(u,\psi)}^*$  is convex-cyclic if and only if  the following holds.
     \begin{enumerate}
\item  $W_{(u,\psi)}^*$
  is cyclic on $\mathcal{F}_2$;
      \item  $|a|=1$, $|u(z_0)|>1$ and $\Im\big(u(z_0) a^m\big)\neq 0$ for all $m\in \NN_0$.
        \end{enumerate}
   \end{lemma}
   \emph{Proof}.  We first assume that  conditions  (i) and (ii) hold,  and proceed to prove the sufficiency.  The condition  $|a|=1$ implies  the operator $W_{(u,\psi)}$ is normal \cite{Tle} and hence $W_{(u,\psi)}^*$ has the same sequence of  eigenvectors $(f_m)$ in \eqref{egv} as $W_{(u,\psi)}$ with corresponding distinct eigenvalues $\overline{u(z_0) a^m}$. Then, following the same argument
   as  in the proof of Lemma~\ref{lem22} and  applying  \cite[Theorem~6.2]{TAN},   the operator is convex-cyclic, and  has a dense set of convex-cyclic vectors.

   Conversely, assume $W_{(u,\psi)}^*$ is convex-cyclic. Then it  is obviously  cyclic. Moreover, the operator $W_{(u,\psi)}$, which is the adjoint of $W_{(u,\psi)}^*$,  has eigenvalues
     $ u(z_0) a^m$  as referred above.  Thus, by \cite[Proposition~3.3]{HR}, it follows that $u({z_0})a^m\in \CC\backslash \big(\overline{\D}\cup \R\big)$ for all $m\in \NN_0$ from which the remaining necessity conditions follow and completes the proof.

If  $u=1$, then    $W_{(u,\psi)}$ is just  the composition  map $C_\psi: f\mapsto f(\psi)$. On the other hand,  if $\psi$ is the identity map, then  $W_{(u,\psi)}$ reduces to the multiplication operator $M_u:  f\mapsto  u\cdot f$. Thus,  $W_{(u,\psi)}$ generalizes the two operators and can be   written as a  product $ W_{(u,\psi)}= M_u C_\psi$.
     The following statement is an immediate consequence of Theorem~\ref{thm22} about  the two factor  operators in the product.
   \begin{corollary}\label{cor1}
\begin{enumerate}
 \item
  Let $C_\psi$ be a bounded composition operator on $\mathcal{F}_2$. Then neither  $C_\psi$ nor its adjoint $C_\psi^*$ is  convex-cyclic on $\mathcal{F}_2$.
  \item Let $M_u$ be a bounded  multiplication operator on $\mathcal{F}_2$. Then $M_u$ can not be convex-cyclic on $\mathcal{F}_2$.
  \end{enumerate}
\end{corollary}
 The proof for the composition operator  follows immediately from Theorem~\ref{thm22} since  $u(z_0)\neq 1$. On the other hand, as proved in \cite[Lemma 2]{CMS}, for $\psi= az+b$, the adjoint of $C_\psi$ is  the  weighted composition operator $C_\psi^*= W_{(K_b,\phi)}$ where $\phi(z)= \overline{a}z$. Then $z_0= 0$ in this case and $|u(z_0)|= |K_b(0)|= 1$. Then the claim follows from condition (ii) of  Theorem~\ref{thm22} again.

By Theorem~\ref{thm1}, the multiplication operator  $M_u$ is not cyclic and can not be convex-cyclic either.

We note in passing that there exists an interesting  interplay between  $u$ and $\psi$ such that $W_{(u,\psi)}= M_u C_\psi$ is bounded (compact) on $\mathcal{F}_2$ while both the factors
    $C_{\psi}$  and $u$ fail to be.  For example one can set  $u_0(z)= e^{-z}$, $\psi_0(z)= z+1$, and observe that   $W_{(u_0,\psi_0)}$ is bounded while both the factors  remain unbounded.   Now Theorem~\ref{thm22} and Corollary~\ref{cor1} provide  another interplay between  $u$ and $\psi$ for which the weighted composition can be convex-cyclic while  both the factors fail to be.
   \subsection{Invariant convex sets for weighted composition operators and their adjoints} \ \  \\
Let  $T$ be a bounded operator on a Banach space $\mathcal{H}$ and  $M$ be a subset of  $\mathcal{H}$.  We say  $M$ is
 invariant  under  $T$ if $T (M)\subseteq M$. We now  study   when the invariant closed convex-sets of the weighted composition operators and their adjoints  are  all invariant subspaces. As shown below, this happens if and only  if the operators  are  convex-cyclic.
  \begin{lemma}\label{lem44}
 Let $u$ and $\psi$  be analytic maps on $\CC$  such that $W_{(u,\psi)}$  is bounded on  $\mathcal{F}_2$, and hence $\psi(z)=az+b, \ |a|\leq 1$.   Then
 \begin{enumerate}
 \item $W_{(u,\psi)}$
 has the property that all of its invariant closed convex-sets are invariant subspaces if and only if
     \begin{enumerate}
\item  $W_{(u,\psi)}$
  is cyclic on $\mathcal{F}_2$;
      \item  $|a|=1$, $|u(z_0)|>1$ and $\Im\big(u(z_0) a^m\big)\neq 0$ for all $m\in \NN_0$.
        \end{enumerate}
        \item $W_{(u,\psi)}^*$
 has the property that all of its invariant closed convex-sets are invariant subspaces if and only if
       $W_{(u,\psi)}$
  is cyclic on $\mathcal{F}_2$ and  $|a|=1$, $|u(z_0)|>1$ and $\Im\big(u(z_0) a^m\big)\neq 0$ for all $m\in \NN_0$.
        \end{enumerate}
        \end{lemma}
       \emph{ Proof}. (i) By  \cite[Proposition 8.11]{NFP}, the operator   $W_{(u,\psi)}$
 has the property that all of its invariant closed convex-sets are invariant subspaces if and only if  for every  closed invariant subspace $M$ of $W_{(u,\psi)}$, the operator
 $W_{(u,\psi)}|M$ is convex-cyclic and the convex-cyclic vectors for $W_{(u,\psi)}|M$ are the same as the cyclic vectors for
$W_{(u,\psi)}|M$ whenever  it is cyclic. Here, by  $W_{(u,\psi)}|M$  we mean the operator  obtained by restricting   $W_{(u,\psi)}$ to  a closed  subset $ M$ of the space
$\mathcal{F}_2$.  In particular, if we   set  $M= \mathcal{F}_2$ and  apply  Lemma~\ref{lem22} above, we observe that  the  necessary conditions in (a) and (b) follow.

 Conversely, assume now that conditions (a) and (b) are satisfied.  Let $M$ be  a closed invariant subspace for  $W_{(u,\psi)}$ and $W_{(u,\psi)}|M $ is cyclic.  We need to show that
$W_{(u,\psi)}|M $ is convex-cyclic.  But this follows readily from Lemma~\ref{lem22} since all eigenvalues of $W_{(u,\psi)}|M $  are  also eigenvalues of $W_{(u,\psi)} $  over the whole space $\mathcal{F}_2$.

The proof of part (ii) follows from a similar argument as  part (i) and Lemma~\ref{lem33}.
  \subsection{Weak and $\tau_{pt}$-supercyclic weighted composition operators  } \ \ \ \\
  Having completely identified the cyclic and convex-cyclic weighted  composition operators and knowing that the space supports no such supercyclic operators, it is natural
to seek weakening the topology of the space and study the supercyclic  structure  with respect to the weak topology (weak supercyclicity ) and the pointwise convergence topology $\tau_{pt}$ ($\tau_{pt}$-supercyclic ).  The  weak and  $\tau_{pt}$-supercyclicities are defined by  simply replacing the norm topology by  these respective topologies on the space. Clearly, weak supercyclicity  is a  stronger property than  $\tau_{pt}$-supercyclicity. The next  diagram exhibit the relations among the various forms of cyclicities for bounded operators. \\
\\

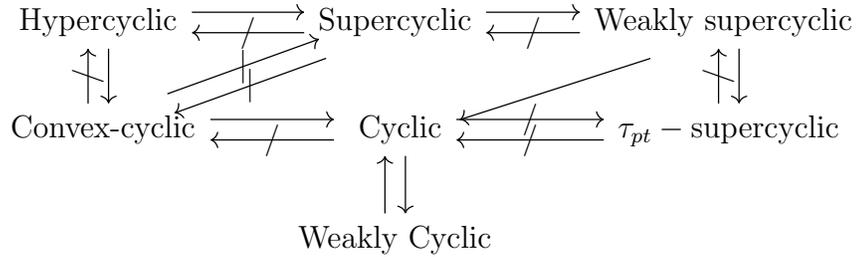
  \begin{tikzcd}
\text{Hypercyclic}  \arrow[r, shift left=.75ex] \arrow[d,shift left=.75ex]
 & \text{Supercyclic} \arrow[r, shift left=.75ex]  \arrow[l,"/" marking, shift left=.75ex] \arrow[dl,"/" marking, shift left=.75ex]
 & \text{Weakly  supercyclic }  \arrow[l, shift left=.75ex, "/" marking] \arrow[d,shift left=.75ex]  \arrow[dl,shift left=.75ex] \\
   \text{ Convex-cyclic}  \arrow[r,  shift left=.75ex] \arrow[u,"/" marking, shift left=.75ex] \arrow[ur,"/" marking, shift left=.75ex]
       & \text{ Cyclic}\arrow[r,"/" marking, shift left=.75ex]      \arrow[l,"/" marking, shift left=.75ex] \arrow[d, shift left=.75ex]
       &\tau_{pt}-\text{supercyclic} \arrow[u, "/" marking, shift left=.75ex] \arrow[l, shift left=.75ex, "/" marking] \\
   & \text{Weakly  Cyclic}  \arrow[u, shift left=.75ex]\\
      \end{tikzcd}
\\

     We note that there exist other forms of  weaker hypercyclicity and  supercyclicity which are not listed in the diagram above. For example  between  supercyclicity and weak supercyclicity, one can find weak l-sequentially supercyclic and weak sequentially  supercyclic properties which the first  implies the second.  A weighted composition operator never satisfies any  of these forms as it fails to satisfy the weakest form of supercyclicity with respect to the pointwise convergence topology as will be seen in Theorem~\ref{supercyclic}. We also note that  since the weak closure of
the convex set, span Orb(T, f) for any $f$ in the given space,  coincides with its norm closure, cyclicity in the norm
topology is equivalent to cyclicity in the weak topology.

Let us now consider the weak  and $\tau_{pt}$-supercyclicity of  the weighted composition operator.  A  simple observation in this regard is that if  $W_{(u,\psi)}$ is weakly supercyclic on $\mathcal{F}_2$, then by \cite[Theorem C]{TAN}, the whole point spectrum of $W_{(u,\psi)}$ and its adjoint must belong to the open disc  $(0, \|W_{(u,\psi)}\|)$. On the other hand, if the operator is convex-cyclic, by Theorem~\ref{thm22},  $\psi(z)= az+b, |a|=1$, then \eqref{simplenorm} implies $ \|W_{(u,\psi)}\|=  |u(0)|e^{|b|^2/2}=|u(z_0)|$. This together with the eigenvalues
gives the necessary condition that
    \begin{align*}
    |a^m u(z_0)|= |a|^m ||u(z_0)|=|u(z_0)| <\|W_{(u,\psi)}\|=|u(z_0)|
\end{align*} for all $ m\in \NN_0$ which is a contradiction. This   assures  that the concepts of  weak supercyclicity and convex-cyclicity  are not related either.  In fact, we will prove that $W_{(u,\psi)}$  is never weakly supercyclic on $\mathcal{F}_2$. The next  result shows that the Fock space supports no $\tau_{pt}$-supercyclic weighted composition operators either.
\begin{theorem}\label{supercyclic}
 Let  $(u,\psi)$ be a pair  of entire functions on $\CC$ which induce a bounded  weighted composition operator $W_{(u,\psi)}$ on $\mathcal{F}_2$. Then  $W_{(u,\psi)}$ can not be  supercyclic on $\mathcal{F}_2$ with respect to  the pointwise convergence topology.
\end{theorem}
As illustrated in the  diagram, pointwise topology is weaker than the weak topology on $\mathcal{F}_2$ and hence the space supports no weakly  supercyclic weighted composition operators. It is interesting to remark that  convex-cyclicity is the strongest form of cyclicity  for weighted composition operators supported on  $\mathcal{F}_2$.

\textbf{Proof of Theorem~\ref{supercyclic}.}
   Since $W_{(u,\psi)}$ is bounded, we set $\psi(z)= az+b$, with $|a|\leq 1$ and consider two different cases.\\
   \emph{Case 1}: Let $|a|<1$ or $|a|=1$ and $a\neq1$  or $a=1$ and $b=0$. For this case the  map $\psi$ fixes the point $z_0= \frac{b}{1-a}$ for $a\neq 1$ and $z_0=0$ for the rest.   Assume on the contrary that there exists a $\tau_{pt}$-supercyclic vector
$f$ in $\mathcal{F}_2$.  First we claim that $u$  is  zero free on $\CC$ because if $u$ vanishes at point $w$, then  \eqref{interplay} implies that every element in the projective orbit of $f$ vanishes at $w$ which extends to the closure resulting a contradiction.  Observe also  that   $f$ can not have zero in $\CC$. if not,  all the elements in the projective orbit will also vanish at a possible zero  which extends to the closure and contradicts.  Thus,  by Proposition~4 of \cite{Belt}, for any two different numbers $z, w\in \CC$,
\begin{align}
\label{dense}
\overline{\bigg\{\frac{u_n(z) f(\psi^n(z))}{u_n(w) f(\psi^n(w))}\bigg\}}= \CC.
\end{align}
Let $r>0$ be given. Then the set   $K= \{z\in \CC: |z-z_0|\leq r\}$ is a compact neighbourhood of $z_0$  which also contains $\psi(K)$ since  for each  $z\in K$
\begin{align*}
|\psi(z)-z_0|= |az+b-z_0| \leq |az-az_0|+|az_0-z_0+b| \leq |a|r\leq r.
\end{align*}
 Now, if we set $w= \psi(z), z\in K, z\neq z_0$ and consider the expression in \eqref{dense}
 \begin{align*}
 \bigg|\frac{u_n(z) f(\psi^n(z))}{u_n(w) f(\psi^n(w))}\bigg|=   \bigg|\frac{u(z) f(\psi^n(z))}{u(\psi^n(z)) f(\psi^{n+1}(z))}\bigg| \leq M
 \end{align*}  for all $n\in \NN$ where
 \begin{align*}
 M= \frac{\max_{z\in K}|u(z)|\cdot\max_{z\in K}|f(z)|}{\min_{z\in K}|u(z)|\cdot \min_{z\in K}|f(z)|}.
 \end{align*} This obviously  contradicts the relation in  \eqref{dense}.

\emph{Case 2}. It remains to show the case for $a=1$ and $b\neq0$.   To this end, let $Q_b= \{z\in \CC: |z-b|\leq 2|b|\}$. Then  $Q_b$  is  a compact  neighbourhood of $b$, and also contains $\psi(Q_b)$ since  for each  $z\in K_b$
  \begin{align*}
|\psi(z)-b|= |z+b-b| \leq |z-b|+|b| \leq 2|b|.
\end{align*}
Then we arrive at the desired conclusion by simply  replacing the compact set $K$ by $Q_b$  in the above argument and completes the proof.

The following corollary is now an immediate consequence of Theorem~\ref{supercyclic}.
  \begin{corollary}
\begin{enumerate}
 \item
  Let $C_\psi$ be a bounded composition operator on $\mathcal{F}_2$. Then neither  $C_\psi$ nor  its adjoint $C_\psi^*$ can  be  $\tau_{pt}$-supercyclic on $\mathcal{F}_2$.
  \item Let $M_u$ be a bounded  multiplication operator on $\mathcal{F}_2$. Then $M_u$ can not be $\tau_{pt}$-supercyclic on $\mathcal{F}_2$.
  \end{enumerate}
  \end{corollary}
Setting $u=1$ in Theorem~\ref{thm1}, we note  that $C_\psi$ is cyclic on $\mathcal{F}_2$  if and only if $a^k\neq a$ for all positive integers $k>1$. The same conclusion can be also read in \cite{KK,TMW}.

\noindent
We end this section with an important remark. We note that the  above proof does not use  Hilbert space property from  $\mathcal{F}_2$. Thus, the same result with the same proof is valid on other classical Fock spaces $\mathcal{F}_p, 1\leq p<\infty $ which consists of all analytic functions $f$  for which
   \begin{align*}
 \|f\|_{ p}^p= \frac{p}{2\pi} \int_{\CC} |f(z)|^p
e^{-\frac{p}{2} |z|^2}  dA(z) <\infty
\end{align*} where $dA$ is the usual  Lebesgue area measure on the complex plane $\CC$. This further  answers  a question raised by   T. Carrol and C. Gilmore in \cite{CC}.
  
\end{document}